\documentclass{amsart}
\usepackage{amssymb,latexsym,graphics,amsfonts}

\swapnumbers
 \theoremstyle{plain}

 \theoremstyle{definition}
 
 \theoremstyle{remark}

 \numberwithin{equation}{subsection}

\newcommand{\f}{\frac}

\newcommand{\ti}{\times}

\newcommand{\op}{\bigoplus}
\newcommand{\ovs}{\overset}
\newcommand{\subs}{\subseteq}
\newcommand{\al}{\alpha}
\newcommand{\del}{\delta}

\newcommand{\lo}{\longrightarrow}

\newcommand{\BN}{\Bbb N}

\newcommand{\BR}{\Bbb R}

 \newcommand{\cal}[1]{\mathcal{#1}}

\begin{document}

\title{Derivations into  duals of closed ideals of Banach algebras }
\author{M. Eshaghi Gordji }
\address{Department of Mathematics,
Semnan University, Semnan, Iran} \email{maj\_ess@Yahoo.com}
\author{B. Hayati }
\address{Department Of Mathematics, Shahid Beheshti University, Tehran, Iran}
\email{b\_hayati@cc.sbu.ac.ir}
\author{S. A. R. Hosseiniun }
\address{Department Of Mathematics, Shahid Beheshti University, Tehran, Iran}
\email{ahossinioun@yahoo.com} 


\dedicatory{}



\smallskip

\begin{abstract}
Let $\cal A$ be a Banach algebra. We study those closed ideals $I$
of $\cal A$ for which the first cohomology group of $\cal A$ with
coefficients in $I^*$ is trivial; i.e.  $H^1(\cal A,I^*)=\{0\}$. We
investigate such closed ideals  when $\cal A$ is weakly amenable or
biflat. Also we give some hereditary properties of ideal
amenability.

\end{abstract}

\maketitle


\section{Introduction}






Let $\cal A$  be a Banach algebra and $X$ a Banach $\cal
A-$bimodule. Then  $X^*$, the dual space of $X$, is also a Banach
$\cal A-$bimodule with module multiplications defined by
\begin{center}
$\langle x,a.x^* \rangle = \langle x.a,x^* \rangle ,
\hspace{0.5cm} \langle x,x^*.a \rangle=\langle a.x,x^* \rangle ,
\hspace{1cm}(a\in {\cal A}, ~x\in X,~x^*\in X^*).$
\end{center}
In particular, $I$ and $I^*$ are  Banach $\cal A-$bimodule  for
every closed ideal $I$ of $\cal A$. A derivation from $\cal A$
into $X$ is a continuous linear operator $D$ such that
$$D(ab)=a\cdot D(b)+D(a)\cdot b \hspace{1cm}(a,b\in {\cal A}).$$

 We define $\delta _{x}(a)=a\cdot x-x\cdot a $ ;  for each $x\in X$ and $a\in \cal A$. $\delta _{x}$ is a derivation from
$\cal A$ into $X$, which is called inner derivation. A Banach
algebra $\cal A$  is amenable if every derivation from $\cal A$
into every dual $\cal A-$bimodule  $X^*$ is inner i.e. $H^1({\cal
A},X^*)=\{0\}$. This definition was introduced by B. E. Johnson in
[J1], [Run1] and [He]. A Banach algebra $\cal A$ is weakly
amenable if every derivation from $\cal A$  into ${\cal A}^*$  is
inner i.e. $H^1({\cal A},{\cal A}^*)=\{0\}$. Bade, Curtis and
Dales [B-C-D] have introduced the concept of weak amenability for
commutative Banach algebras. You can see also [J3], [D-Gh], [G1],
[G2] and [G3]. Let $n\in \mathbb N$, a Banach algebra $\cal A$ is
called n-weakly amenable if $H^1({\cal A},{\cal A}^{(n)})=\{0\}$,
where ${\cal A}^{(n)}$ is the n-th dual of $\cal A$ [D-Gh-G] and
[J2].

Let $G=SL(2,\BR)$, the set of elements in ${\Bbb M}_2( \BR)$ with
determinant 1 , also let ${\cal A}=L^1(G)$ and $I$ be the
augmentation ideal of ${\cal A}$, then theorem 5.2 of [J-W] implies
that $H^1({\cal A},I^*)\neq\{0\}$. On the other hand ${\cal A}$ is
weakly amenable. This example guides us to the following
definitions.

Let $\cal A$ be a Banach algebra and $I$ be a closed two-sided ideal
of $\cal A,$ then $\cal A$ is $I-$weakly amenable if every
derivation from $\cal A$ into $I^*$ is inner, in other words
$H^1({\cal A},I^*)=\{0\}$.
 We call $\cal A$  ideally amenable if $\cal A$
is $I-$weakly amenable for every closed ideal $I$ of $\cal A$ [E-Y]
, [E-H]. Let $n\in \Bbb N$, a Banach algebra $\cal A$ is called
n-ideally amenable if for every closed two-sided ideal $I$ in $\cal
A$, $H^1({\cal A},I^{(n)})=\{0\}.$

Obviously, an ideally amenable Banach algebra is weakly amenable.
Since every closed ideal of $\cal A$ is a Banach $\cal A-$bimodule,
then  an amenable Banach algebra is ideally amenable. There are some
examples of Banach algebras to show that ideal amenability is not
equivalent to weak amenability or amenability. In the following we
give some of them.

1- Let $\cal A$ be the unitization of the augmentation ideal of
$L^1(SL(2,\BR)).$ Then $\cal A$ is weakly amenable and $\cal A$ is
not ideally amenable [E-Y].

2- Let $\cal A=L^1(SL(2,\BR)).$ Then $\cal A$ is weakly amenable and
$\cal A$ is not ideally amenable.

3- Let $\cal A$ be a non-nuclear $C^*$-algebra. Then $\cal A$ is
non-amenable, ideally amenable Banach algebra [E-Y].

4- Let $\cal A$ be a commutative non-amenable, weakly amenable
Banach algebra. Then $\cal A$ is non-amenable, n-ideally amenable
Banach algebra for each $n\in \Bbb N$ [E-Y].

 Let $\cal A$ be a Banach algebra , $X$  a Banach  $\cal A-$bimodule and $Y$  a
closed $\cal A-$submodule of $X$, we say that the short exact
sequence $ \{0\} \lo Y\ovs{i}{\lo} X \ovs{\pi}{\lo} \f{X}{Y} \lo
\{0\}$ of $\cal A-$bimodules  splits if $\pi$ has a bounded right
inverse which is also an $\cal A-$bimodule homomorphism. The
following theorem is well known.
\paragraph{\bf Theorem 1.1}Let $\cal A$ be a Banach algebra, $X$
a Banach $\cal A-$bimodule and  $Y$  a closed $\cal A-$submodule of
$X$. Then the following conditions are equivalent.
\begin{itemize}
\item[i)]The short exact sequence
$\{0\} \lo Y\ovs{i}{\lo}  X \ovs{\pi}{\lo} \f{X}{Y} \lo \{0\}$
splits.
\item[ii)]$i$ has a bounded left inverse which is also an $\cal A-$bimodule
homomorphism.
\item[iii)]There exists a continuous projection of $X$ onto $Y$
which is also an $\cal A-$bimodule homomorphism.
\end{itemize}
See [D] for a proof.

 Let $\cal A$ be a Banach algebra. Then the projective tensor
product of $\cal A$ is denoted by $\cal A\hat{\otimes}_\pi\cal A$
. This space is a Banach $\cal A$-bimodule with module
multiplications defined by :
\begin{center}
 $ a.(b\otimes c)=ab\otimes c \hspace{1cm}$ and $ \hspace{1cm}(b\otimes c).a=b\otimes ca\hspace{2cm} (a,b,c\in \cal A)$.\\
\end{center}
 The corresponding diagonal operator $\bigtriangleup :\cal A\hat{\otimes}_\pi\cal A \lo \cal A $ is defined
by $a\otimes b\mapsto ab$ . It is clear that $\bigtriangleup$ is a
Banach $\cal A$-bimodule homomorphism.

  Let $E$ be a Banach space and $F(E)$ be the space of finite
rank operators on $E$. We say that $E$ has the approximation
property if there is a net $(S_\alpha)_\alpha$ in $F(E)$ such that
$(S_\alpha)\lo id_E$ uniformly on compact subsets of $E$.
\section{derivations into duals of submodules}
 Let $X$ be a Banach  ${\cal A}-$bimodule, and  $Y$  a closed ${\cal
A}-$submodule of $X$. By using exact sequences, we give some
conditions that $H^1({\cal A},X^*)=\{0\}$ implies  $H^1({\cal
A},Y^*)=\{0\}.$
\paragraph{\bf Theorem 2.1} Let ${\cal A}$ be a Banach algebra,
$X$ a Banach  ${\cal A}-$bimodule and  $Y$  a closed ${\cal
A}-$submodule of $X$. If  $H^1({\cal A},X^*)=\{0\}$ and the exact
sequence
\begin{align*}
\{0\} \lo Y^\perp \ovs{i}{\lo} X^* \ovs{\pi}{\lo} \f{X^*}{Y^\perp}
\lo \{0\}\hspace{2cm}(1)
\end{align*}
of Banach $\cal A-$bimodules splits,  then $H^1({\cal
A},Y^*)=\{0\}$.
\paragraph{\bf Proof.} Let $D:{\cal A}\lo Y^*$ be a derivation. Since the exact sequence
 (1) splits, $\pi$ has a bounded right inverse, say $\phi$, such that $\phi$ is also an
${\cal A}-$bimodule homomorphism. In this case $\phi\circ D:{\cal
A}\lo X^*$ is a derivation, so there exists $f\in X^*$ such that
$\phi\circ D=\del_f$. Thus, we have $\pi\circ \phi\circ D=\pi\circ
\del_f$. This shows that $id_{Y^*}\circ D=\del_{\pi(f)}$ and
therefore $D=\del_{\pi(f)}$.\hfill$\blacksquare~$\\

 The following lemma is in literature but we give its proof.
\paragraph{\bf Lemma 2.2.}
The exact sequence
\begin{align*}
\{0\} \lo Y^\perp \ovs{i}{\lo} X^* \ovs{\pi}{\lo} \f{X^*}{Y^\perp}
\lo \{0\}\hspace{2cm}(1)
\end{align*}
 splits, if the following exact sequence splits ;
\begin{align*}
\{0\}\lo Y\ovs{i}{\lo} X\ovs{\pi}{\lo} \f{X}{Y}\lo \{0\}\hspace
{2cm}(2)
\end{align*}
\paragraph{\bf Proof.} Since the exact sequence (2) splits, there exists a continuous  projection
$P$ of $X$ onto $Y$ which is also an ${\cal A}-$bimodule
homomorphism. Let $Q=id_{X^*}-P^*$. Then for each $y\in Y$ and
$f\in X^*$ we have
\begin{align*}
\langle y,Q(f)\rangle  &=\langle y,f-P^*f\rangle  =\langle y,f\rangle -\langle y,P^*f\rangle  \\
&= \langle y,f\rangle -\langle P(y),f\rangle =0.
\end{align*}
 So $Q(X^*)\subseteq Y^\perp$. On the other hand for each  $f\in X^*$ and $x\in X$ we have
\begin{align*}
\langle x,Q^2(f)\rangle  &=\langle x,Q(f-P^*f)\rangle = \langle x,(f-P^*f-P^*(f-P^*f))\rangle  \\
&=\langle x,f\rangle -\langle P(x),f\rangle -\langle P(x),f\rangle  +\langle P(x),P^*f\rangle \\
&=\langle x,f\rangle -\langle P(x),f\rangle = \langle x,Q(f)\rangle.
\end{align*}
Thus, $Q$ is a continuous projection of $X^*$ onto $Y^\perp$. Also
for $f\in X^*$, $a\in {\cal A}$ and $x\in X$, we have
\begin{align*}
\langle x,Q(a.f)\rangle  &=\langle x,a.f-P^*(a.f)\rangle = \langle x.a,f\rangle  -\langle P(x),a.f\rangle  \\
&=\langle x.a,f\rangle  -\langle P(x).a,f\rangle =\langle x.a,f\rangle  -\langle P(x.a),f\rangle \\
&=\langle x.a,Q(f)\rangle =\langle x,a.Q(f)\rangle .
\end{align*}
 So $Q$ is a left ${\cal A}-$module  homomorphism. Similarly
$Q$ is a right $\cal A-$module homomorphism and this completes the
proof .\hfill$\blacksquare~$
\paragraph{\bf Corollary 2.3.} Let $\cal A$, $X$, $Y$ be as in Theorem 2.1 . If the exact sequence (2) splits and
 $H^1({\cal A},X^*)=\{0\}$ , then $H^1({\cal A},Y^*)=\{0\}$.
\paragraph{\bf Corollary 2.4.} Let ${\cal A}$ be a Banach algebra and $n\in\BN$. If
 $H^1({\cal A},X^{(n+2)})=\{0\}$, then  $H^1({\cal
A},X^{(n)})=\{0\}$.
\paragraph{\bf Proof.} Let $\wedge_{n-1}:X^{(n-1)}\lo X^{(n+1)}$ be the canonical map.
Then the exact sequence
$$\{0\}\lo X^{(n-1)}\ovs{\wedge _{n-1}}{\lo} X^{(n+1)}
\ovs{\pi}{\lo} \f{X^{(n+1)}}{X^{(n-1)}} \lo \{0\}$$ splits, because
the adjoint of  $\wedge _{n-2}$ ,   $ \wedge _{n-2}^*:X^{(n+1)}\lo
X^{(n-1)} $, is a left inverse of $\wedge
_{n-1}$ which is also an $\cal A-$bimodule homomorphism. Now use corollary 2.3 .\hfill$\blacksquare~$\\

 The next corollary has been proved in [D-Gh-G], but it is an immediate result of Corollary 2.4.
\paragraph{\bf Corollary 2.5.} Let ${\cal A}$ be $n+2-$ weakly amenable ($n+2-$ideally amenable). Then ${\cal A}$ is
$n-$weakly amenable ($n-$ideally amenable).
\section{ Closed ideals of weakly amenable Banach algebras }
In this section, we find some closed ideals of a weakly amenable
Banach algebra $\cal A$ for which  $H^1(\cal
 A,I^*)$ is trivial. We denote the linear span of the set
  $\lbrace ab:a,b\in \cal A\rbrace$ by ${\cal A}^2$.
  We  show  that  if a closed ideal $I$
 satisfies ${\cal A}^2\subseteq I$ and  $H^1(\cal A,I^*)=\{0\}$,
 then ${\cal A}^2$ is dense in $I$. This is a generalization of
 Gr\o nb\ae k's theorem [D, Theorem 2.8.63].

 First by  Theorem 2.1 and Corollary 2.3 we have the following
 theorem.
\paragraph{\bf Theorem 3.1} Assume that ${\cal A}$ is a weakly amenable Banach
algebra. If  one of the following conditions holds for each closed
ideal $I$ in ${\cal A}$, then  ${\cal A}$ is ideally amenable .
\begin{itemize}
\item[i)] The exact sequence $ \{0\} \lo I^\perp \ovs{i}{\lo} \cal
A^* \ovs{\pi}{\lo} \f{\cal A^*}{I^\perp} \lo \{0\}$, splits.
\item[ii)] The exact sequence $\{0\}\lo I\ovs{i}{\lo} \cal
A\ovs{\pi}{\lo} \f{\cal A}{I}\lo \{0\}$, splits.
\end{itemize}
\paragraph{\bf Theorem 3.2.}Let ${\cal A}$ be a Banach algebra
and  $I$ be a closed ideal of ${\cal A}$ with a bounded
approximate identity. Then the following conditions are equivalent
\begin{itemize}
\item[i)] $I$ is weakly amenable;
\item[ii)]${\cal A}$ is $I-$weakly amenable.
\end{itemize}
\paragraph{\bf Proof.} Let ${\cal A}$ be $I-$weakly amenable and
$D:I\rightarrow I^*$  a derivation. Since $I$ is psudo unital $\cal
A-$bimodule, by Proposition 2.1.6 of [Run1] $D$  has an extension
$\bar{D}:{\cal A}\rightarrow I^*$ such that $\bar{D}$ is a
derivation . But $\cal A$ is $I-$weakly amenable, thus $\bar{D}$ and
consequently $D$ is inner. The converse is Lemma 2.1 of
[E-Y].\hfill$\blacksquare~$\\

We recall that, in a Banach algebra $\cal A$, a net  $(e_\al)_\al$
is quasi-central if for each element $a\in \cal A$ ;
$\displaystyle\lim _ \al (ae_\al-e_\al a) =0$  . Obviously, each
approximate identity is a quasi-central net.
\paragraph{\bf Theorem 3.3.} Let ${\cal A}$ be a weakly amenable
Banach algebra and  $I$  a closed ideal of ${\cal A}$ with a
quasi-central bounded  approximate identity . Then ${\cal A}$ is
$I-$weakly amenable.
\paragraph{\bf Proof.} Let $(e_\al)$ be a quasi-central bounded  approximate identity in
$I$ and let $J$ be an ultrafilter on the index set of $(e_\al)$
such that dominates the order filter. Define
\begin{center}
$P: {\cal A}^*\lo {\cal A}^*$ \\
$\phi~{\longmapsto} w^*-\displaystyle\lim_J  (\phi-e_\al.\phi)$\\
\end{center}
 For every $\phi\in {\cal A}^*$ and $a\in I$ we have
$$\langle a,P\phi\rangle = \lim_J \langle a,\phi-e_\al.\phi\rangle
=\lim_J \langle a-ae_\al,\phi\rangle =0.$$ Thus $PA^*\subs I^\perp$.
Also for $\phi\in I^\perp$ and $a\in {\cal A}$, we have
$$\langle a,P\phi\rangle =\lim_J (\langle a,\phi\rangle -\langle ae_\al,\phi\rangle )
=\langle a,\phi\rangle .$$
This means that $P$ is a projection of ${\cal A}^*$ onto
$I^\perp$. On the other hand for $a,b\in {\cal A}$ and $\phi\in
{\cal A}^*$ we have
\begin{align*}
\langle b,P(a.\phi)\rangle  &=\lim_J \langle b,a.\phi-e_\al.(a.\phi)\rangle  \\
&=\lim_J \langle  ba-be_\al a,\phi \rangle =\lim_J \langle ba-bae_\al,\phi\rangle \\
&=\lim_J \langle ba,\phi-e_\al.\phi\rangle =\langle ba,P(\phi)\rangle  \\
&=\langle b,a.P(\phi)\rangle .
\end{align*}
So, $P$ is a left ${\cal A}-$module  homomorphism. Similarly $P$
is a right ${\cal A}-$module homomorphism. Therefore the exact
sequence $\{0\}\lo I^\perp\ovs{i}{\lo} {\cal A}^*\ovs{\pi}{\lo}
\f{{\cal A}^*}{I^\perp} \lo \{0\}$ splits and consequently ${\cal
A}$ is $I-$weakly amenable.\hfill$\blacksquare~$

Let $\cal A$ be a Banach algebra and let $I$ be a closed  ideal of
$\cal A$. If $I$  is Arens regular with bounded approximate identity
$(e_\alpha)_{\al}$, then $(e_\alpha)_{\al}$
 is quasi-central. To prove this, let $E\in I^{**}$ be a cluster
 point of $(e_\alpha)_{\al}$. Then $E$ is the identity of
 $I^{**}$. Also for each $a\in \cal A,$ we have $aE,Ea\in
 I^{**}$ and thus $aE=E(aE)=(Ea)E=Ea.$

So, by applying the above theorem, we have the following.

\paragraph{\bf Corollary 3.4.} Let ${\cal A}$ be a weakly amenable Banach
algebra and let $I$ be a closed ideal of $\cal A$ with a bounded
approximate identity.  If one of the following conditions holds,
\begin{itemize}
\item[i)] $I$ is Arens regular.
\item[ii)]$\cal A$ is Arens
regular.
\end{itemize}
Then ${\cal A}$ is $I-$weakly amenable.

Since every $C^*$-algebra $\cal A$ is weakly amenable [Ha], Arens
regular, and every closed ideal of $\cal A$ has a bounded
approximate identity, then for each closed ideal $I$ of $\cal A,$
we have $H^1(\cal A, I^*)=\{0\}.$ In the other words, every
$C^*$-algebra is ideally amenable [E-Y].
\paragraph{\bf Corollary 3.5.} Let ${\cal A}$ be a weakly amenable Banach algebra
 such that  each closed ideal of ${\cal A}$ has a quasi-central bounded
  approximate identity. Then ${\cal A}$ is ideally amenable.\\

 Let $(a_\al)_{\al\in I}$ be a quasi-central bounded net in $\cal A$. Then the
following closed ideals of ${\cal A}$ are called the net ideals of
${\cal A}$  [E].
\begin{quote}
$\displaystyle I(a_\alpha):=\{a\in {\cal A}~;~~\lim_\alpha aa_\alpha= a~\}$,\\
$\displaystyle K(a_\alpha):=\{a\in {\cal A}~;~~\lim_\alpha aa_\alpha=0~\}$,\\
$C(a_\alpha):=\{a\in {\cal A}~;~~(aa_\alpha)\ is~~ convergent \}$,\\
$\displaystyle L(a_\alpha):=\{\lim_\alpha aa_\alpha~;~~a\in
C(a_\alpha)~\}.$
\end{quote}
\paragraph{\bf Theorem.3.6.} Let $I(a_\al),K(a_\al),C(a_\al)$ and $L(a_\al)$ are as above.
 If $L(a_\al)=I(a_\al)$ and ${\cal A}$ is
$C(a_\al)-$weakly amenable, then
\begin{itemize}
\item[i)] ${\cal A}$ is $I(a_\al)-$weakly amenable. \item[ii)]
${\cal A}$ is $K(a_\al)-$weakly amenable.
\end{itemize}
\paragraph{\bf Proof.} Since $L(a_\al)=I(a_\al)$, we have $C(a_\al)=I(a_\al)\op K(a_\al)$  [E]. Consequently
 the exact sequences $\{0\}\lo {I(a_\al)}\lo C(a_\al)\lo \f{C(a_\al)}{I(a_\al)} \lo \{0\}$
and $\{0\}\lo {K(a_\al)}\lo {C(a_\al)}\lo \f{C(a_\al)}{K(a_\al)}
\lo \{0\}$
 split. Now by Corollary 2.3, ${\cal A}$ is both $I(a_\al)-$weakly amenable and  $K(a_\al)-$weakly amenable.\hfill$\blacksquare~$
 \paragraph{\bf Theorem.3.7.} Let $(a_\al)_{\al\in I}$ be a quasi-central bounded net in $\cal A$. If one of the
following assertions holds, then ${\cal A}$ is $I(a_\al)-$weakly
amenable.
\begin{itemize}
\item[i)] ${\cal A}$ is weakly amenable and $\{a_\al: \al\in
I\}\subs I(a_\al)$.  \item[ii)] There exists a codimension one
ideal $M$ of ${\cal A}$ such that $H^1(M,I(a_\al)^*)=\{0\}$.
\end{itemize}
\paragraph{\bf Proof.} If (i) holds, then $(a_\al)_{\al\in I}$ is a quasi-central
bounded right approximate identity  for $I(a_\al)$. By Theorem
3.3, ${\cal A}$ is $I(a_\al)-$weakly amenable. Let (ii) holds.
Since $(a_\al)_{\al\in I}$ is a bounded approximate identity in
${\cal A}$ for ${\cal A}$-bimodule $I(a_\al)$, then by Cohen's
factorization theorem, we have $AI(a_\al)=I(a_\al){\cal
A}=I(a_\al)$, so by [G-L, Theorem 2.3], $H^1(\cal A
,I(a_\al)^*)=H^1(M,I(a_\al)^*)$. Thus ${\cal A}$ is
$I(a_\al)-$weakly amenable.\hfill$\blacksquare~$\\

We know that if $\cal A$ is a weakly amenable Banach algebra, then
$\bar{{\cal A}^2}=\cal A$ [D, Theorem 2.8.63]. A generalization of
this fact is as follows.
\paragraph{\bf Theorem 3.8.} Let ${\cal A}$ be a Banach algebra and  ${\cal A}^2\subseteq
I$ for an arbitrary closed two-sided ideal $I$ of ${\cal A}$. If
${\cal A}$ is $I-$weakly amenable, then ${\cal A}^2$ is dense in
$I$.
\paragraph{\bf Proof.} Assume $\bar{{\cal A}}^2\neq I$. Then there exists
$0\neq \varphi\in I^*$ such that $\varphi|_{\bar{{\cal A}^2}}=0$.
Let $\Phi\in {\cal A}^*$ be a Hahn-Banach extension of $\varphi$
on ${\cal A}$. We define $D:{\cal A}\lo I^*$ by $a\mapsto \langle
a,\Phi\rangle \varphi$. For $a,b\in {\cal A}$, $D(ab)=\langle
ab,\Phi\rangle \varphi=0$. Also for $i\in I$ we have $\langle
i,D(a).b\rangle =\langle i,(\langle a,\Phi\rangle
\varphi).b\rangle =\langle a,\Phi\rangle \langle bi,\varphi
\rangle =0$. So $D(a).b=0$ . Similarly $a.D(b)=0$. Therefore $D$
is a derivation from ${\cal A}$ into $I^*$ and by hypothesis there
exists $\psi\in I^*$ such that $D=\del_\psi$. Also for every $i\in
I$ we have
\begin{align*}
\langle i,\varphi\rangle ^2&=\langle i,\varphi\rangle  \langle
i,\Phi\rangle =\langle i,D(i)\rangle \\
&=\langle i,\del_\psi(i)\rangle  =\langle i,i\psi-\psi i\rangle \\
&=0.
\end{align*}
 So $\varphi=0$, which is  contradiction. Thus ${\cal A}^2$ is
 dense in $I$.\hfill$\blacksquare~$
\paragraph{\bf Corollary 3.9.} Let $\cal A$ be a Banach algebra
and let $M$ be a closed non maximal modular ideal of $\cal A$ with
codimension one. If $H^1(\cal A,M^*)=\{0\}$, then ${\cal A}^2$ is
dense in $M$.
\paragraph{\bf Proof.} Since the codimension of $M$ is one, there
exists $a\in \cal A$ such that $\cal A=\Bbb {C}a+M$. We show that
$a^2 \in M$. Assume that $a^2$ does not belong to $M$, so there
exists $ 0\neq\alpha$ and $m\in M$ such that $a^2=\alpha a+m$. Now
let $b$ be an arbitrary element of $\cal A$, there exist $\beta$
and $m'\in M$ such that $b=\beta a+m'$, so
\begin{align*}
b-b({\alpha}^{-1}a)&=(\beta a+m')-(\beta a+m')({\alpha}^{-1}a)\\
&=m'-\beta{\alpha}^{-1}m-{\alpha}^{-1}m'a.
\end{align*}
Thus, for each $b$, $b-b({\alpha}^{-1}a)$ belongs to $M$. This
means that ${\alpha}^{-1}a$ is a left modular identity for $M$.
Similarly ${\alpha}^{-1}a$ is a right modular identity for $M$, so
$M$ is a maximal modular ideal which is contradiction. Therefore
$a^2$ belongs to $M$ and consequently $\cal A^2\subseteq M$. Now
by the above theorem, ${\cal A}^2$ is dense in
$M$.\hfill$\blacksquare~$

\section{Closed ideals of biflat Banach algebras}

We say that a Banach algebra $\cal A$ is biprojctive if
$\bigtriangleup :\cal A \hat{\otimes}_\pi\cal A\lo \cal A$ has a
bounded right inverse which is an $\cal A$-bimodule homomorphism.
Also we say that a Banach algebra $\cal A$ is biflat if the bounded
linear map ${\bigtriangleup}^*:{\cal A}^*:\lo ({\cal A
\hat{\otimes}_\pi\cal A})^*$ has a bounded left inverse which is an
$\cal A$-bimodule homomorphism [Run1]. Obviously  by taking
adjoints, one sees that every biprojective Banach algebra is biflat.
It is well known that every biflat Banach algebra is weakly amenable
[D], and a Banach algebra is amenable if and only
if it is biflat  and has a bounded approximate identity [Run2]. \\

 Since there is no Hahn-Banach theorem for operators, there is none
for bilinear continuous forms. In other words, let $E$ and $F$ be
two Banach spaces and $G$ be a subspace of $E$ and $\phi \in
BL(G,F;\Bbb C)$, where $BL(G,F;\Bbb C)$ is the set of all bounded
bilinear mappings from $G\times F$ into $\Bbb C$. In general case,
there is no any extension of $\phi$ to a bilinear map $\tilde
{\phi}\in BL(E,F;\Bbb C)$. Since $BL(G,F;\Bbb C)\simeq L(G,F^*)$
this situation is equivalent to say that each element $T\in
L(G,F^*)$ doesn't have any extension to an element $\tilde {T}\in
L(E,F^*)$ [D-F,1.5].

 However, there is some conditions that Hahn-Banach theorem
works for operators as well.\\
Let  $\pi(z;E,F)$ be the projective norm of the element  $z\in
F\hat{\otimes}_\pi F$ and $G$ be a subspace of $E$. Then it is
clear that $\pi (z;E,F)\leq \pi (z;G,F)$ for each element $z\in
G\hat{\otimes}_\pi F$ . If there exists $\lambda \geq 1$ such that
$\pi (z;G,F)\leq \lambda \pi (z;E,F)$ for each element $z\in
G\hat{ \otimes}_\pi F$, then we say that $.\hat{\otimes}_\pi F$
respects $G$ into $E \hat{\otimes}_\pi F$ isomorphically . For
example, $.\hat{\otimes}_\pi F$ respects $G$ into $E
\hat{\otimes}_\pi F$ isomorphically if $G$ is a complemented
subspace of $E$ [D-F].

 By Hahn-Banach theorem we can extend each element of $T\in(G\hat{\otimes}_\pi F)^*$
 to a continuous linear functional $\tilde{T}$ on $E \hat{\otimes}_\pi F$ provided
that $.\hat{\otimes}_\pi F$ respects $G$ into $E \hat{\otimes}_\pi
F$ isomorphically [D-F].

Now, we are ready to bring our main theorem about biflat Banach
algebras. Before doing this we recall that an ideal $I$ is left
essential as a left Banach $\cal A$-module if the linear span of
$\lbrace ai : ~~~a\in \cal A,~~~i\in I \rbrace$ is dense in $I$.\\
\paragraph{\bf Theorem.4.1.} Let $\cal A$ be a biflat Banach
algebra and $I$ a closed ideal of $\cal A$ which is left essential.
If
  $\cal A \hat{\otimes}_\pi .$   respects $I$ into $\cal A
\hat{\otimes}_\pi \cal A$ isomorphically, then $H^1(\cal
A,I^*)=\{0\}$.
\paragraph{\bf Proof.}Let $D:\cal A \lo I^*$ be a derivation.
Since $\cal A$ is biflat, ${\bigtriangleup}^*:{\cal A}^*\lo ({\cal
A \hat{\otimes}_\pi\cal A})^*$ has a bounded left inverse $\rho$
which is an $\cal A$-bimodule homomorphism. Let
\begin{center}
$T:L(\cal A,I^*)\lo (\cal A \hat{\otimes}_\pi I)^*$\\
 $S\longmapsto T_S$
\end{center}
 be the isometric isomorphism which is defined by $\langle a\otimes
 i,T_S \rangle=\langle i,S(a)\rangle$. Let $\widetilde {T}_D$ be a
 Hahn-Banach extension of $T_D$ on $\cal A\hat {\otimes}_\pi \cal
 A$ . We claim that $D=\delta_\phi$ , where $\phi =\rho
 (\widetilde{T}_D)|_I$ .\\
 First we show that for each $i\in I$ and $a\in \cal A$ ;
$i.(a.\widetilde {T}_D - \widetilde {T}_D.a)=i.\bigtriangleup ^*
(\widetilde {Da})$,
 where $\widetilde {Da}$ is a Hahn-Banach extension of
$Da$ on $\cal A$. Let $b,c\in \cal A$, we have
\begin{align*}
\langle b\otimes c,i.(a.\widetilde {T}_D - \widetilde
{T}_D.a)\rangle &=\langle b\otimes ci,a.\widetilde {T}_D -
\widetilde {T}_D.a\rangle \\
&=\langle b\otimes cia-ab\otimes ci,\widetilde{T}_D\rangle \\
&=\langle b\otimes cia,T_D\rangle-\langle ab\otimes ci,T_D\rangle \\
&=\langle cia,Db\rangle - \langle ci,D(ab)\rangle \\
&=\langle ci,a.Db-D(ab)\rangle =\langle ci,Da.b\rangle \\
&=\langle bci,Da\rangle =\langle bci,\widetilde{Da}\rangle \\
&=\langle b\otimes ci,\bigtriangleup ^*(\widetilde {Da})\rangle \\
&=\langle b\otimes c,i.\bigtriangleup ^*(\widetilde {Da})\rangle .
\end{align*}
Now, let $a,b\in \cal A$ and $i\in I$. Then we have
\begin{align*}
\langle bi,\delta_\phi (a)\rangle &=\langle bi,a.\phi-\phi .a\rangle\\
&=\langle bi,a.\rho(\widetilde{T}_D)|_I-\rho(\widetilde{T}_D)|_I.a \rangle \\
&=\langle bia-abi,\rho(\widetilde{T}_D)|_I \rangle \\
&=\langle bia-abi,\rho(\widetilde{T}_D)\rangle \\
&=\langle bi,a.\rho(\widetilde{T}_D)-\rho(\widetilde{T}_D).a\rangle \\
&=\langle b,\rho(i.(a.\widetilde{T}_D-\widetilde{T}_D.a))\rangle \\
&=\langle b,\rho(i.\bigtriangleup ^* (\widetilde {Da}))\rangle \\
&=\langle bi,id_{{\cal A}^*} (\widetilde {Da})\rangle .
\end{align*}
Since $I$ is left essential and $Da$ , $\delta_\phi(a)$ are both
continuous linear functional on $I$, we have $Da=\delta_\phi(a)$.
This is true for each $a\in\cal A$, so $D=\delta_\phi$ and $D$ is
inner.\hfill$\blacksquare~$
\paragraph{\bf Corollary.4.2.}Let $\cal A$ be a biflat Banach
algebra with a left approximate identity. Then $\cal A$ is ideally
amenable provided that  $\cal A\hat{\otimes}_\pi .$  respect all
closed ideals into $\cal A\hat{\otimes}_\pi\cal A$ isomorphically.

There are a kind of biprojective Banach algebras whose left closed
ideals are left essential. These algebras are semiprime biprojctive
Banach algebras with the approximation property [S].
\paragraph{\bf Corollary.4.3.}Let $\cal A$ be a semiprime,
biprojective Banach algebra with the approximation property, and
$I$ a closed ideal of $\cal A$. If  $\cal A\hat{\otimes}_\pi .$
 respects $I$ into $\cal A\hat{\otimes}_\pi\cal A$ isomorphically,
then $H^1(\cal A,I^*)=\{0\}$. In particular, for each closed ideal
$I$ which is complemented as a subspace of $\cal A$, the assertion
holds.

\section{Some Hereditary properties of ideal amenability}

Let $\cal A$ and $\cal B$ be two Banach algebras and $\phi:\cal
A\lo \cal B$ a continuous homomorphism with dense range. We know
that $\cal B$ is amenable if $\cal A$ is amenable [J1], but this
is not true for weak amenability. In special case, if $\cal A$ is
weakly amenable and commutative, then  $\cal B$ is weakly amenable
[D].

\paragraph{\bf Theorem 5.1} Let $\cal A$ and $\cal B$ are two Banach
algebras and let   $\phi:\cal A\lo \cal B$ be a continuous
homomorphism with dense range and  $J$ a closed ideal of $\cal B$.
If the following conditions hold:
\begin {itemize}
\item[i)]$\phi |_{J^c}$ is one to one, where $J^c$ is
$\phi^{-1}(J)$. \item[ii)]$\phi(J^c)$ is dense in $J$.
\item[iii)]$H^1(\cal A,{J^c}^*)=\{0\}$.
\end{itemize}
 Then $H^1(\cal B,J^*)=\{0\}$.
\paragraph{\bf Proof.} Let $D:\cal B\lo J^*$  be a derivation . Define $T:J^*\lo {J^c}^*$ , given by $f\mapsto T_f$,
where $T_f$ is defined by $T_f(a)=f(\phi(a))$. Obviously  for each
$f\in J^*$ , $T_f$ is a continuous linear functional on $J^c$ and
so $T$ is well defined. We show that $T$ is onto.\\
Let $g\in {J^c}^*$ and define $\tilde {f}:\phi(J^c)\lo \Bbb {C}$
by $\phi(a)\mapsto g(a)$. $\phi(J^c)$ is a subspace of $J$ and
$\tilde{f}$ is a bounded linear functional on $\phi(J^c)$. By
Hahn-Banach theorem, there exists a continuous linear functional
$f$ on $J$ such that $f|_{\phi(J^c)}=\tilde{f}$. It is clear that
$T_f=g$.\\
Now define $\tilde{D}:\cal A\lo {J^c}^*$ by $\tilde{D}:=T\circ
D\circ \phi$. $\tilde{D}$ is a bounded linear map. Let $a_1,a_2\in
\cal A$ and $a\in J^c$. Then
\begin{align*}
\langle \tilde{D}(a_1a_2),a\rangle &=\langle
T(D(\phi(a_1a_2))),a\rangle \\
&=\langle
T(D(\phi(a_1)).\phi(a_2)+\phi(a_1).D(\phi(a_2))),a\rangle\\
&=\langle D(\phi(a_1)).\phi(a_2)+\phi(a_1).D(\phi(a_2)),\phi(a)\rangle\\
&=\langle D(\phi(a_1)),\phi(a_2a)\rangle+\langle D(\phi(a_2)),\phi(aa_1)\rangle\\
&=\langle T\circ D\circ \phi(a_1),a_2a\rangle +\langle T\circ D\circ
\phi(a_2),aa_1\rangle\\
&=\langle \tilde {D}(a_1).a_2,a\rangle +\langle
a_1.\tilde{D}(a_2),a\rangle
\end{align*}
Thus $\tilde{D}$ is a derivation. Since $H^1(\cal
A,{J^c}^*)=\{0\}$ , then there exists $g\in {J^c}^*$ such that
$\tilde{D}=\delta_g$. But T was onto so, there exists $f\in J^*$
such that $T_f=g$. We claim that $D=\delta_f$ . Let $a\in \cal A$
and $ a_1\in J^c$, We have
\begin{align*}
\langle D(\phi(a)),\phi(a_1)\rangle&=\langle
\tilde{D}(a),a_1\rangle=\langle a.g-g.a,a_1\rangle\\
&=\langle g,a_1a-aa_1\rangle\\
&=\langle T_f,a_1a-aa_1\rangle\\
&=\langle f,\phi(a_1a-aa_1)\rangle\\
&=\langle\phi(a).f-f.\phi(a),\phi(a_1)\rangle\\
&=\langle \delta_f(\phi(a)),\phi(a_1)\rangle
\end{align*}
Since $\phi(J^c)$ is dense in $J$, we have
$D(\phi(a))=\delta_f(\phi(a))$ for each $a\in \cal A$. Again since
$\phi(\cal A)$ is dense in $\cal B$, then
$D=\delta_f$.\hfill$\blacksquare~$

Let  $\cal A$ and $\cal B$ be two Banach algebras and $\phi:\cal
A\lo \cal B$ a continuous homomorphism with dense range. In general,
 we assert that the ideal amenability of $\cal A$ doesn't imply the
ideal amenability of $\cal B$.

We know that, the approximation property is not necessary for the
weak amenability of the algebra of approximable operators on a
Banach space [B,Corollary.3.5]. Also there are some Banach spaces
$E$ with the approximation property such that $A(E)$ is not weakly
amenable [B,Theorem.5.3].

Now let $E$ be a Banach space with the approximation property such
that $A(E)$ is not weakly amenable. Then the nuclear algebra
$N(E)$ of $E$ is biprojective and consequently weakly amenable.
Since $N(E)$ is topologically simple, then $N(E)$ is ideally
amenable. On the other hand, $A(E)$ is not ideally amenable and
the inclusion map $i:N(E)\lo A(E)$ is a continuous homomorphism
with dense range. This proves the assertion.

\paragraph{\bf Theorem.5.2.} Suppose  $Y$ and $Z$ are closed
subspaces of a Banach space $X$, and suppose that there is a
collection $\Lambda \subset B(X)$ with the following properties:
\begin{itemize}
\item[i)] Every $\phi\in \Lambda$ maps $X$ into $Y$.
\item[ii)] Every $\phi\in \Lambda$ maps $Z$ into $Z$.
\item[iii)]$\sup \lbrace \parallel \phi\parallel: \phi\in
\Lambda\rbrace < \infty$.
\item[iv)] To every $y\in Y$ and to every $\epsilon >0$ corresponds a
$\phi\in \Lambda$ such that $\parallel y-\phi y\parallel <\epsilon$.
\end{itemize}
Then $Y+Z$ is closed.
\paragraph{\bf Proof.}[Rud, 1.2.Theorem].

\paragraph{\bf Corollary.5.3.}
Suppose $\cal A$ is a Banach algebra. Let $I$ be a right closed
ideal and $J$ be a left closed ideal of $\cal A$. If $I$ has a
bounded approximate identity, then $I+J$ is closed.
\paragraph{\bf Proof.} Let $(e_\al)_\al$ be a bounded approximate identity for $I$ and
$$ \Lambda=\lbrace L_{e_\al}:\cal A\lo \cal A\hspace{.1cm}  |\hspace{0.5cm} L_{e_\al}(a)=e_\al
a\hspace{.3cm}\rbrace .$$  Obviously $\Lambda\subset B(\cal A)$. For
each $\al$, we have $L_{e_\al}(\cal A)=\lbrace e_\al a: a\in \cal
A\rbrace\subset I$  and  $L_{e_\al}(J)=\lbrace e_\al j:j\in
J\rbrace\subset J$. Also we have  $$\sup \lbrace \parallel
L_{e_\al}\parallel : L_{e\al}\in \Lambda\rbrace \leq\displaystyle
\sup_\al\lbrace \parallel e_\al\parallel \rbrace <\infty.$$ Now let
$\epsilon >0$ is given and $i\in I$. There exists $\al_0$ such that
$$\parallel i-L_{e_{\al_0}}(i)\parallel= \parallel i-e_{\al_0}i\parallel
 <\epsilon.$$
Thus by the above theorem, $I+J$ is closed.\hfill$\blacksquare~$
\paragraph{\bf Theorem.5.4.} Let $\cal A$ be a Banach algebra and
let $I$ be a closed ideal of $\cal A$ with bounded approximate
identity. If $I$ and $\f {\cal A}{I}$ are ideally amenable, then
$\cal A$ is ideally amenable.
\paragraph{\bf Proof.} Let $J$ be a closed ideal of $\cal A$ and $D:
\cal A\lo J^*$ a derivation. Consider $\imath :I\cap J\lo J$ as an
inclusion map.Obviously $\imath^*:J^*\lo (I\cap J)^*$ is an $\cal
A$-bimodule homomorphism and so, $\imath^*\circ D:\cal A\lo(I\cap
J)^*$ and consequently  $\imath^*\circ D |_I:I\lo(I\cap J)^*$ is a
continuous derivation. Since $I$ is ideally amenable, there exists
$\phi_1 \in (I\cap J)^*$ such that $\imath^*\circ
D|_I=\delta_{\phi_1}$. Let $\Phi_1$ be the Hahn-Banach extension
of $\phi_1$ on $J$. Define $\tilde{D}:=D-\delta_{\Phi_1}$ . So,
$\tilde{D}$ is a derivation from $\cal A$ into $J^*$ . We show
that $\tilde{D}|_I=0$. Let $i\in I$ and $j\in J$ , we have
\begin{align*}
\langle j,\tilde{D}(i)\rangle &=\langle
j,D(i)\rangle-\langle j,\delta_{\Phi_1}(i)\rangle \\
&=\langle j,D(i)\rangle -\langle ji-ij,\Phi_1\rangle\\
&=\langle j,D(i)\rangle-\langle ji-ij,\phi_1\rangle\\
&=\langle j,D(i)\rangle-\langle j,\delta_{\phi_1}(i)\rangle\\
&=\langle j,D(i)\rangle-\langle j,\imath^*(D(i))\rangle\\
&=\langle j,D(i)\rangle-\langle \imath(j),D(i)\rangle=0.
\hspace{2cm}
\end{align*}
So, $\tilde{D}|_I=0$ and $\tilde{D}$ induces a map from $\f{\cal
A}{I}$ into $J^*$, we call it $\tilde{D}$ itself, which is a
derivation. Since $I\subset Ann (\f{J}{J\cap I})$, the annihilator
of $(\f{J}{J\cap I})$, then
$\f{J}{J\cap I}$ is a Banach $\f{\cal A}{I}$-bimodule.\\
On the other hand , let $(e_\al)_\al$  be  a bounded approximate
identity in $I$. Then for each $x\in I\cap J$ and $a\in\cal A$ we
have
\begin{align*}
\langle x,\tilde{D}(a+I)\rangle &=\displaystyle\lim
_\al\langle xe_\al,\tilde{D}(a)\rangle\\
&=\displaystyle \lim_\al\langle e_\al,\tilde{D}(a).x\rangle\\
&=\displaystyle \lim_\al \langle e_\al,\tilde{D}(ax)-a.\tilde{D}(x)\rangle\\
&=0.\hspace{2cm}
\end{align*}
So, $\tilde{D}(\f{\cal A}{I})\subseteq (I\cap
J)^\perp\cong(\f{J}{J\cap I})^*$. Since $I$ has a bounded
approximate identity, then by Corollary 5.3 , $J+I$ is a closed
ideal of $\cal A$, thus $\f{J+I}{I}$ is a Banach space. Now define
$\psi:\f{J}{J\cap I}\lo\f{J+I}{I}$ by $j+J\cap I\mapsto j+I$.
Obviously $\psi$ is an algebra isomorphism. Therefore $\psi$ is an
$\f{\cal A}{I}$-bimodule isomorphism. Also
\begin{align*}
\parallel \psi(j+J\cap I)\parallel &=\parallel j+I\parallel\\
&=\inf \lbrace \parallel j+i\parallel:i\in I\rbrace\\
&\leq\inf\lbrace\parallel j+i\parallel :i\in I\cap J\rbrace \\
&=\parallel j+J\cap I\parallel .
\end{align*}
So, $\psi$ is bounded. By open mapping theorem $\psi$ is a
homeomorphism and consequently $(\f{J}{J\cap I})^*\cong
(\f{J+I}{I})^*$. Therefore there exists $\Phi_2\in(I\cap J)^\perp $
such that $\tilde{D}=\delta_{\Phi_2}$ . It shows that
$D=\delta_{\Phi_1+\Phi_2}$ ; $D$ is inner and $\cal A$ is ideally
amenable.\hfill$\blacksquare~$

It is notable that the above theorem is true for weak amenability
and amenability  even if $I$ does not have any bounded approximate
identity [D].

 Now we pose an open problem in this direction.

 \paragraph {\bf Question.} Is valid the above theorem, if $I$ does not
 have any bounded approximate identity ?






\end{document}